\documentclass[11pt]{amsart}
\usepackage{amssymb,amsmath,amsgen,amsxtra,amsfonts,a4wide} 
\usepackage{dsfont,times}

\DeclareFontEncoding{OT2}{}{} 

\begin{document}
%
%
%
\theoremstyle{definition}
\newtheorem{Definition}{Definition}[section]
\newtheorem*{Definitionx}{Definition}
\newtheorem{Convention}{Definition}[section]
\newtheorem{Construction}{Construction}[section]
\newtheorem{Example}[Definition]{Example}
\newtheorem{Examples}[Definition]{Examples}
\newtheorem{Remark}[Definition]{Remark}
\newtheorem{Remarks}[Definition]{Remarks}
\newtheorem{Caution}[Definition]{Caution}
\newtheorem{Conjecture}[Definition]{Conjecture}
\newtheorem{Question}[Definition]{Question}
\newtheorem{Questions}[Definition]{Questions}
\theoremstyle{plain}
\newtheorem{Theorem}[Definition]{Theorem}
\newtheorem*{Theoremx}{Theorem}
\newtheorem{Proposition}[Definition]{Proposition}
\newtheorem*{Propositionx}{Proposition}
\newtheorem{Lemma}[Definition]{Lemma}
\newtheorem{Corollary}[Definition]{Corollary}
\newtheorem{Fact}[Definition]{Fact}
\newtheorem{Facts}[Definition]{Facts}
\newtheoremstyle{voiditstyle}{3pt}{3pt}{\itshape}{\parindent}%
{\bfseries}{.}{ }{\thmnote{#3}}%
\theoremstyle{voiditstyle}
\newtheorem*{VoidItalic}{}
\newtheoremstyle{voidromstyle}{3pt}{3pt}{\rm}{\parindent}%
{\bfseries}{.}{ }{\thmnote{#3}}%
\theoremstyle{voidromstyle}
\newtheorem*{VoidRoman}{}

%
\newcommand{\prf}{\par\noindent{\sc Proof.}\quad}
\newcommand{\blowup}{\rule[-3mm]{0mm}{0mm}}
\newcommand{\cal}{\mathcal}
\newcommand{\Aff}{{\mathds{A}}}
\newcommand{\BB}{{\mathds{B}}}
\newcommand{\CC}{{\mathds{C}}}
\newcommand{\FF}{{\mathds{F}}}
\newcommand{\GG}{{\mathds{G}}}
\newcommand{\HH}{{\mathds{H}}}
\newcommand{\NN}{{\mathds{N}}}
\newcommand{\ZZ}{{\mathds{Z}}}
\newcommand{\PP}{{\mathds{P}}}
\newcommand{\QQ}{{\mathds{Q}}}
\newcommand{\RR}{{\mathds{R}}}
\newcommand{\Liea}{{\mathfrak a}}
\newcommand{\Lieb}{{\mathfrak b}}
\newcommand{\Lieg}{{\mathfrak g}}
\newcommand{\Liem}{{\mathfrak m}}
\newcommand{\ideala}{{\mathfrak a}}
\newcommand{\idealb}{{\mathfrak b}}
\newcommand{\idealg}{{\mathfrak g}}
\newcommand{\idealm}{{\mathfrak m}}
\newcommand{\idealp}{{\mathfrak p}}
\newcommand{\idealq}{{\mathfrak q}}
\newcommand{\idealI}{{\cal I}}
\newcommand{\lin}{\sim}
\newcommand{\num}{\equiv}
\newcommand{\dual}{\ast}
\newcommand{\iso}{\cong}
\newcommand{\homeo}{\approx}
\newcommand{\mm}{{\mathfrak m}}
\newcommand{\pp}{{\mathfrak p}}
\newcommand{\qq}{{\mathfrak q}}
\newcommand{\rr}{{\mathfrak r}}
\newcommand{\pP}{{\mathfrak P}}
\newcommand{\qQ}{{\mathfrak Q}}
\newcommand{\rR}{{\mathfrak R}}
%
%
\newcommand{\dq}{{``}}
\newcommand{\OO}{{\cal O}}
\newcommand{\into}{{\hookrightarrow}}
\newcommand{\onto}{{\twoheadrightarrow}}
\newcommand{\Spec}{{\rm Spec}\:}
\newcommand{\BigSpec}{{\rm\bf Spec}\:}
\newcommand{\Proj}{{\rm Proj}\:}
\newcommand{\Pic}{{\rm Pic }}
\newcommand{\Br}{{\rm Br}}
\newcommand{\NS}{{\rm NS}}
\newcommand{\chit}{\chi_{\rm top}}
\newcommand{\KanDiv}{{\cal K}}
\newcommand{\perdef}{{\stackrel{{\rm def}}{=}}}
\newcommand{\Cycl}[1]{{\ZZ/{#1}\ZZ}}
\newcommand{\Sym}{{\mathfrak S}}
\newcommand{\Xcan}{X_{{\rm can}}}
\newcommand{\Ycan}{Y_{{\rm can}}}
\newcommand{\ab}{{\rm ab}}
\newcommand{\Aut}{{\rm Aut}}
\newcommand{\Hom}{{\rm Hom}}
\newcommand{\Supp}{{\rm Supp}}
\newcommand{\ord}{{\rm ord}}
\newcommand{\divisor}{{\rm div}}
\newcommand{\Alb}{{\rm Alb}}
\newcommand{\Jac}{{\rm Jac}}
\newcommand{\piet}{{\pi_1^{\rm \acute{e}t}}}
\newcommand{\Het}[1]{{H_{\rm \acute{e}t}^{{#1}}}}
\newcommand{\Hcris}[1]{{H_{\rm cris}^{{#1}}}}
\newcommand{\HdR}[1]{{H_{\rm dR}^{{#1}}}}
\newcommand{\hdR}[1]{{h_{\rm dR}^{{#1}}}}
\newcommand{\Elf}{{\rm Elf}}
\newcommand{\defin}[1]{{\bf #1}}
%
\newcommand{\textcyr}[1]{%
 {\fontencoding{OT2}\fontfamily{wncyr}\fontseries{m}\fontshape{n}\selectfont #1}}
\newcommand{\Sha}{{\mbox{\textcyr{Sh}}}}
%

\title[Non-reduced Picard schemes]{A note on non-reduced Picard schemes}
\author{Christian Liedtke}
\address{Mathematisches Institut, Heinrich-Heine-Universit\"at, 40225
  D\"usseldorf, Germany}
\email{liedtke@math.uni-duesseldorf.de}
\thanks{2000 {\em Mathematics Subject Classification.} 14K30, 14J10, 14C20} 
\date{May 7, 2008, {\it revised:} August 26, 2008}

\begin{abstract}
  The Picard scheme of a smooth curve and a smooth complex variety is reduced.
  In this note we discuss which classes of surfaces in terms of the
  Enriques--Kodaira classification can have non-reduced Picard schemes and
  whether there are restrictions on the characteristic of the ground field. 
  It turns out that non-reduced Picard schemes are uncommon in 
  Kodaira dimension $\kappa\leq0$, that this phenomenon can be bounded for
  $\kappa=2$ (general type) and that it is as bad as can be in $\kappa=1$.
\end{abstract}
\setcounter{tocdepth}{1}
\maketitle
\tableofcontents
\section*{Introduction}

The set of isomorphism classes of invertible sheaves on a scheme
$X$ forms a group, the so-called Picard group $\Pic(X)$.
In case $X$ is integral and projective over a field $k$, 
this group $\Pic(X)$ carries a natural scheme structure
as was shown by Grothendieck \cite{gr}.
Moreover, if $X$ is geometrically normal, then
$\Pic^0(X)$, the identity component of $\Pic(X)$, is even
projective.

A theorem of Cartier states that group schemes over fields of
characteristic zero are reduced.
It follows that $\Pic^0$ of a projective and geometrically
normal scheme is an Abelian variety in this case.

However, over fields of positive characteristic, the $\Pic^0$
even of a smooth projective variety need no longer be
reduced.
A first example has been constructed by Igusa \cite{ig}.
As explained by Mumford in \cite[Lecture 27]{mum}, the 
non-reducedness of the Picard scheme can be related to 
Bockstein operations in cohomology. 
It follows that varieties with $h^2(X,\OO_X)=0$ have
a reduced Picard scheme.
And in particular, $\Pic^0$ of a geometrically normal curve 
is always an Abelian variety.

Hence we have to look at dimension two and in view of the
Enriques--Kodaira classification it is natural to ask:
\begin{enumerate}
 \item What kind of surfaces, e.g. ruled, elliptic, or general type, 
   have a non-reduced $\Pic^0$ ?
 \item Fixing numerical invariants,
    is it true that surfaces with these invariants have a reduced $\Pic^0$?  
  \item If the previous question has a negative answer in general, does it have
    a positive answer if the characteristic of the ground field is
    sufficiently large?
\end{enumerate}
%

From the Enriques--Kodaira classification and its extension to
positive characteristic by Bombieri--Mumford \cite{bm2} we
immediately get

\begin{Propositionx}
  For Kodaira dimension $\kappa(X)\leq0$, the Picard scheme tends to be reduced:
  \begin{enumerate}
    \item If $\kappa(X)=-\infty$ then $\Pic^0(X)$ is reduced.
    \item If $\kappa(X)=0$ then $\Pic^0(X)$ is reduced 
     except for a few exceptional cases
     in characteristic $2$ and $3$.
   \end{enumerate} 
\end{Propositionx}

In Kodaira dimension $\kappa=1$ all surfaces possess
elliptic fibrations and the
non-reducedness of the Picard scheme is closely related to
the existence of wild fibres.
Using results on torsors under Jacobian fibrations we show the
following, which is more or less implicit in the literature:

\begin{Theoremx}
  Let $f:X\to B$ be an elliptic fibration 
  of a surface 
  in positive characteristic.
  Assume that $f$ is not generically constant.
  Then there exists an elliptic fibration $f':X'\to B$ such that
  \begin{enumerate}
   \item $\Pic^0(X')$ not reduced,
   \item $\kappa(X')=1$.
   \item the Jacobian fibrations of $f$ and $f'$ coincide, and
   \item $b_i(X)=b_i(X')$ for all $i$ and $\chi(\OO_X)=\chi(\OO_{X'})$.
  \end{enumerate}
\end{Theoremx}

In particular, for every positive characteristic and
every set of Betti-numbers, Euler characteristic and
not generically constant elliptic fibration for which there 
exists a surface with $\kappa=1$, there exists a surface with the same 
invariants and a non-reduced Picard scheme.

Moreover, we can choose the difference between $h^{01}$ and $\frac{1}{2}b_1$, 
which can be viewed as a measure of the non-reducedness of
$\Pic^0$, as large as we want to.

Examples of Katsura and Ueno show that the situation is similarly bad for
iso-trivial fibrations.
\medskip

For Kodaira dimension $\kappa=2$, i.e., surfaces of general type, there are
examples due to Serre with non-reduced Picard schemes 
in every characteristic.
However, we can limit this phenomenon

\begin{Theoremx}
 Given an integer $m$, there exists an integer $P(m)$, such
 that minimal surfaces of general type with $K_X^2=m$
 over fields of characteristic $p\geq P(m)$
 have a reduced $\Pic^0$.
\end{Theoremx}

\begin{VoidRoman}[Acknowledgements]
 I thank Matthias Sch\"utt for comments as well as the referee for comments 
 and simplifying the proof of Theorem \ref{allgemeinertyp}.
\end{VoidRoman}

\section{Kodaira dimension at most zero}

Let $X$ be a smooth projective surface over an algebraically 
closed field $k$.
We denote by $\kappa(X)$ its Kodaira dimension.
Thanks to the Enriques--Kodaira classification
that was extended to positive characteristic by
Bombieri and Mumford
we have an explicit description of surfaces with
$\kappa(X)\leq0$, which allows us to answer the questions
posed in the introduction quite satisfactory.

Two smooth projective surfaces that are birational are related
by a sequence of blow-ups and blow-downs in closed points.
Since this process does not affect $\Pic^0$, we may and will
restrict ourselves to suitable minimal models in the following.

\begin{Theorem}
 \label{offensichtlich}
 If $\kappa(X)=-\infty$ then $\Pic^0(X)$ is reduced.
\end{Theorem}

\prf
A surface with $\kappa(X)=-\infty$ is birational to $\PP^1\times C$,
where $C$ is a smooth curve.
Hence such a surface has a reduced $\Pic^0$.
\qed


\begin{Theorem}
 \label{kodairanull}
 If $\kappa(X)=0$ then
 $\Pic^0(X)$ is reduced except possibly if $X$ is 
 \begin{enumerate}
  \item a non-classical Enriques surface in characteristic $2$, or
  \item a (quasi-) hyperelliptic surface in characteristic $2$ or $3$.
 \end{enumerate}
 The exceptions do occur. 
\end{Theorem}

\prf
A look at the table of possible invariants 
in the introduction of \cite{bm2} shows that the only surfaces
with non-reduced $\Pic^0$ (noted as $\Delta\neq0$ in this table)
are non-classical Enriques surfaces or certain (quasi-)hyperelliptic 
surfaces.

Non-classical Enriques surfaces can exist in characteristic $2$ only
\cite[Theorem 5]{bm2} and such surfaces have been constructed 
in \cite[Section 3]{bm3}.

Hyperelliptic surfaces of Kodaira dimension zero with non-reduced
$\Pic^0$ are those with $p_g=1$, using the table of possible invariants
again.
These are precisely the hyperelliptic surfaces where $K_X$ is of order 
$1$, and the
detailed analysis in \cite[Section 3]{bm2} shows that such surfaces 
can and do exist in characteristic $2$ and $3$ only.

Quasi-hyperelliptic surfaces exist in characteristic $2$, $3$ only
\cite{bm3}.
As explained in the proof of \cite[Proposition 8]{bm3}, such a surface
has ${\rm ord}K_X=1$, i.e., a non-reduced $\Pic^0$, 
if and only if the character
$K\to\Aut(C_0)/\GG_a\cdot A\iso\GG_m$ is trivial
(notation as in loc.cit.).
In characteristic $2$, this condition is fulfilled for surfaces
of type $f)$ of the {\it Char. $2$}-table in \cite[page 214]{bm3}.
In characteristic $3$, this condition holds for surfaces
of type $d)$ of the {\it Char. $3$}-table in \cite[page 214]{bm3}, 
cf. also \cite[Section 3B]{la}.
\qed

\section{Elliptic fibrations}

We have seen in the first section that surfaces with $\kappa\leq0$ and
non-reduced $\Pic^0$ form a very small class.
This is not true in Kodaira dimension $\kappa=1$, even when
fixing numerical invariants.
Since all these surfaces are endowed with an elliptic 
fibration we translate our
problem into the language of elliptic fibrations.
In fact, twisting an elliptic fibration and adding wild fibres
we can make its $\Pic^0$ as non-reduced as we want to whilst
fixing numerical invariants.

We recall that $H^1(\OO_X)$ can be identified with the Zariski tangent space
to $\Pic^0(X)$ and that $\frac{1}{2}b_1(X)$ is the dimension of $\Pic^0(X)$.
Hence the difference $h^{01}-\frac{1}{2}b_1$ can be viewed as a measure
for the non-reducedness of $\Pic^0$, which is zero if and only if 
$\Pic^0$ is reduced.

Let $f:X\to B$ be an elliptic fibration over a curve.
We recall that a fibre $F$ is called {\em wild}, if 
$h^0(F,\OO_F)\geq2$.
Wild fibres can exist over fields of positive characteristic
only and we refer to 
\cite[Chapter V]{cd} for details.
The following result explains the role of wild fibres 
in view of non-reduced Picard schemes.

\begin{Proposition}
 \label{wild}
 Let $f:X\to B$ be a relatively minimal elliptic fibration over a curve $B$.
 \begin{enumerate}
  \item If $\chi(\OO_X)\geq1$ and $f$ has no wild fibres then $\Pic^0(X)$ is reduced.
  \item If $f$ has $w\geq2$ wild fibres then $\Pic^0(X)$ is not reduced and
   $h^{01}-\frac{1}{2}b_1\geq (w-1)$.
 \end{enumerate}
\end{Proposition}

\prf
We have $R^1f_\ast\OO_X\iso{\cal L}\oplus{\cal T}$, where $\cal L$
is a line bundle on $B$ and $\cal T$ is a torsion sheaf whose support 
consists precisely of those points over which the fibre of $f$ is wild.
From the Grothendieck--Leray spectral sequence we obtain a short exact
sequence
\begin{equation}
 \label{gl}
  0\,\to\,H^1(B,\,\OO_B)\,\to\,H^1(X,\,\OO_X)\,\to\,H^0(B,\,R^1f_\ast\OO_X)\,\to\,0\,.
\end{equation}
Assume $\chi(\OO_X)\geq1$ and that $f$ has no wild fibres.
Then $h^0({\cal T})=0$ and 
the canonical bundle formula for elliptic fibrations gives
$\deg{\cal L}=-\chi(\OO_X)\leq-1$, hence $h^0(B,{\cal L})=0$.
We obtain $h^1(\OO_X)=h^1(\OO_B)$.
By its universal property, the composition $X\to B\to{\rm Jac}(B)$
factors over the Albanese variety of $X$, from which we conclude
$b_1(X)\geq b_1(B)=2h^1(\OO_B)=2h^1(\OO_X)$.
Since we have $b_1(X)\leq2h^1(\OO_X)$ in any case, we obtain
$2h^1(\OO_X)=b_1(X)$, which implies that $\Pic^0(X)$ is
reduced.

Now, assume that $f$ has $w\geq2$ wild fibres.
Then $h^0({\cal T})\geq w$ and hence $h^1(\OO_X)-h^1(\OO_B)\geq w$
by (\ref{gl}).
By \cite[Lemma 3.4]{ku}, we have 
$\frac{1}{2}b_1(X)\leq h^1(\OO_B)+1$, which yields the
desired inequality.
Since $h^{01}$ is strictly larger than $\frac{1}{2}b_1$,
the $\Pic^0(X)$ is non-reduced.
\qed\medskip

The next result tells us that, given an elliptic surface in positive
characteristic that is not generically constant, we can always 
find another fibration with $\kappa=1$ and with the same Betti numbers but 
with arbitrary non-reduced Picard scheme.
In particular, we cannot bound the non-reducedness by fixing
invariants or the characteristic.

\begin{Theorem}
 \label{drehen}
 Let $f:X\to B$ be a relatively minimal elliptic fibration over a 
 curve $B$ defined
 over an algebraically closed field of positive characteristic.
 Let $n\geq1$ be an integer and assume that $f$ is not
 generically constant.

 Then there exists an elliptic fibration $f':X'\to B$, such that
 \begin{enumerate}
  \item $\chi(\OO_X)=\chi(\OO_{X'})$, $K^2_X=K^2_{X'}=0$ and 
      $b_i(X)=b_i(X')$ for all $i$,
  \item both elliptic fibrations have the same Jacobian fibration, 
  \item $\Pic^0(X')$ is non-reduced and even $h^{01}-\frac{1}{2}b_1\geq n$, and
  \item $\kappa(X')=1$ if $n\geq2$ or $p\geq5$.
 \end{enumerate}
\end{Theorem}

\prf
We use the notation of \cite[Section 5.4]{cd}.
We denote by $j:J\to B$ the Jacobian fibration associated to $f:X\to B$.
Let $J_\eta^\sharp$ be the N\'eron model of the generic fibre of $j$.
We denote by $\Elf(j)$ the abelian group classifying torsors under
$J_\eta^\sharp$.
For every  closed point $b\in B$, we let $\tilde{\OO}_{B,b}$ be the
(strict) Henselisation of the local ring $\OO_{B,b}$.
Let $\tilde{J}_b^\sharp$ be the N\'eron model of (the reduction) of
$J\times_B\Spec\tilde{\OO}_{B,b}$ and let
$\Elf(\tilde{j}_b)$ be the abelian group of torsors under 
$\tilde{J}_b^\sharp$.
For every closed point $b\in B$ there exists a homomorphism
$\psi_b:\Elf(j)\to\Elf(\tilde{j}_b)$, the so-called
{\em local invariant}.

Since $f$ is not generically constant, $j$ is not trivial
and in this case there exists a short exact sequence
\begin{equation}
 \label{twisting}
 0\,\to\,\Sha(J_\eta^\sharp)\,\to\,\Elf(j)
 \,\stackrel{\psi}{\to}\,
 \bigoplus_{b\in B}\, \Elf(\tilde{j}_b)
 \,\to\,0\,,\mbox{ \quad where\quad }\psi=\sum_{b\in B}\psi_b\,,
\end{equation}
cf. \cite[Proposition 5.4.3]{cd} and \cite[Corollary 5.4.6]{cd}.

The generic fibre of $j$ is an ordinary elliptic curve as $j$ is
not trivial.
If the fibre above $b$ is an ordinary elliptic curve, there
exists a non-trivial subgroup 
$\Elf(\tilde{j}_b)^{\rm rad}$ of $\Elf(\tilde{j}_b)$, 
such that an element of $\Elf(j)$,
which maps to a non-trivial element of
$\Elf(\tilde{j}_b)^{\rm rad}$ 
corresponds to an elliptic fibration 
with Jacobian fibration $j$ and a wild
fibre above $b$, cf.
\cite[Corollary 5.4.3]{cd}.

We choose a set $S$ of $(n+1)$ distinct points in $B$ such that the fibres
of $j$ above these points are ordinary elliptic curves.
For every $b\in S$ we choose a non-trivial element $e_b$ in 
$\Elf(\tilde{j}_b)^{\rm rad}$.
By the surjectivity of $\psi$ in (\ref{twisting}), there exists
an element $f'$ of $\Elf(j)$ such that $\psi_b(f')=e_b$ 
for every $b\in S$.
This $f'$ corresponds to an elliptic fibration $f':X'\to B$ with
wild fibres above $S$.
By Proposition \ref{wild}, we have $h^{01}-\frac{1}{2}b_1\geq n$ and
that $\Pic^0(X')$ is not reduced.

By \cite[Proposition 5.3.6]{cd} we have $\chi(\OO_X)=\chi(\OO_J)=\chi(\OO_{X'})$ 
and the same for the Betti numbers and $c_2$
by \cite[Corollary 5.3.5]{cd}.
We have $K^2=0$ in any case.

If $h^{01}-\frac{1}{2}b_1\geq n\geq1$ then $\kappa(X')\geq0$ 
by Theorem \ref{offensichtlich}.
By the table of possible invariants in the introduction of \cite{bm2},
we see that $\kappa(X')=0$ and $n\geq1$ implies $h^{01}-\frac{1}{2}b_1=1$ and $p\leq3$.
Hence if $n\geq2$ or $p\geq5$ we have $\kappa(X')=1$.
\qed\medskip

Even among iso-trivial elliptic surfaces with $\kappa=1$ we find
arbitrary non-reduced Picard schemes
in arbitrary positive characteristic.
The following examples are due to Katsura and Ueno:

\begin{Proposition}
 For every prime $p$ and every integer $n$ there exists an 
 elliptic surface with $\kappa=1$ defined over an
 algebraically closed field of characteristic $p$ such that
 \begin{enumerate}
  \item the elliptic fibration is iso-trivial
  \item $\Pic^0$ is not reduced and even $h^{01}-\frac{1}{2}b_1\geq n$.
 \end{enumerate}
\end{Proposition}

\prf
Let $X$ be an elliptic surface of \cite[Example 8.1]{ku}.
As $X$ possesses an iso-trivial elliptic fibration, we 
have $\chi(\OO_X)=0$.
By \cite[Lemma 3.5]{ku} we have $b_1=2$.
For $m\geq3$ (as defined in \cite[Example 8.1]{ku}) 
we have $\kappa(X)=1$ and choosing $m$ sufficiently large,
we get $p_g$ as large as we want to, i.e., we also
get $h^{01}$ as large as we want to since $\chi(\OO_X)=0$.
\qed

\section{General type}

There exist surfaces with $\kappa=2$, i.e., surfaces of general type, 
with non-reduced Picard schemes in arbitrary large characteristic.
However, fixing $K_X^2$, there exists only a finite number of 
characteristics where minimal surfaces of general type with 
these invariants can have non-reduced Picard schemes.

We recall that surfaces of general type can have 
non-reduced Picard schemes in arbitrary large characteristic -
the examples are due to Serre:

\begin{Proposition}
 \label{witt}
 For every prime $p>0$ there exists a minimal surface of general type
 over an algebraically closed field of characteristic $p$ that has
 a non-reduced $\Pic^0$.
\end{Proposition}

\prf
In \cite[Proposition 15]{se}, Serre constructs for every $p>0$ a 
smooth hypersurface $Y_p$ in $\PP^3$ with a fixed point free 
action of $C_p:=\ZZ/p\ZZ$.
From the construction it is clear that we may assume that 
$Y_p$ is of degree $\geq5$, i.e., of general type.
Thus, the quotient $X_p:=Y_p/C_p$ is a surface of general type
with $h^1(\OO_{X_p})\neq0$ by \cite[Proposition 16]{se}.
On the other hand, $b_1(Y_p)=0$ implies $b_1(X_p)=0$ since 
$C_p$ acts without fixed points.
Hence $\Pic^0$ is not reduced.
\qed

\begin{Remark}
 Examples of uniruled surfaces of general type in characteristic
 $2$ with arbitrary non-reduced $\Pic^0$ have been constructed
 in \cite[Theorem 8.1]{uniruled}.
\end{Remark}

\begin{Theorem}
 \label{allgemeinertyp}
 Given an integer $m$, there exists an integer $P(m)$, such
 that minimal surfaces of general type with $K_X^2=m$
 over fields of characteristic $p\geq P(m)$
 have a reduced $\Pic^0$.
\end{Theorem}

\prf
Fixing $K_X^2$, 
the Euler characteristic $\chi(\OO_X)\leq1+p_g$ is
bounded above by Noether's inequality
and bounded below $\chi(\OO_X)\geq0$ 
in characteristic $p\geq11$ by \cite[Theorem 8]{sb}.
Hence there is only a finite number of possibilities for
$\chi(\OO_X)$ if $p\geq11$.

Canonical models of surfaces of general type with 
fixed $\chi(\OO_X)$ and $K_X^2$
are parametrised by a subset of an appropriate Hilbert scheme which 
is defined over $\Spec\ZZ$.
Hence there exists a scheme ${\cal M}$ 
of finite type over $\Spec\ZZ$ and a 
family $f:{\cal X}\to{\cal M}$ of 
canonical models of surfaces of general type such that every 
such surface with $K^2=m$ occurs in this family.

There exists an integer $P_1$ such that for every prime
$p\geq P_1$ all components of ${\cal M}_p$ are flat over 
$\Spec\ZZ$.
Let ${\cal M}'$ be one of these finitely many components.
By \cite{ar} there exists a quasi-finite morphism
${\cal N}'\to{\cal M}'$ and a family 
$f':{\cal Y}\to{\cal N}'$ that resolves the
singularities of $f$ simultaneously.

Then, ${\cal N}'\otimes_\ZZ\QQ$ is non-empty and parametrises 
smooth and minimal surfaces of general type in characteristic zero.
By the Lefschetz principle, we
may assume that the family $f'$ is defined over the
complex numbers.
Then, by Ehresmann's fibration theorem, these surfaces are
diffeomorphic, which implies that all of them have
the same first Betti number $b_1$.
Hence $h^{01}$ is constant in this family being
equal to $b_1/2$ by Hodge theory.
It follows that not only the $\Pic^0$ of all fibres in
this family over ${\cal N}'\otimes_\ZZ\QQ$
are reduced but that also $h^{01}$ is constant.

By upper semicontinuity there exists a closed subset 
${\cal V}\subseteq{\cal N}'$ over which $h^{01}$ of a fibre may
jump.
By Chevalley's theorem, the image of ${\cal V}$ in $\Spec\ZZ$ is 
a constructible set, i.e., closed or open since $\Spec\ZZ$ is
one-dimensional.
However, by what we have just seen, this image avoids the
generic point of $\Spec\ZZ$ and so this image is a proper
closed subset.
In particular, there exists a $P_2'$, such that
for every prime $p\geq P_2'$, the fibre ${\cal N}_p'$ does
not intersect with ${\cal V}$.
Since $p\geq P_1$, for every field 
$K$ of characteristic $p\geq{\rm max}(P_1,P_2')$ and every
morphism $\Spec K\to{\cal N}'$ the fibre 
${\cal Y}_K:={\cal Y}\times_{{\cal N}'}\Spec K$ is a 
surface of general type that lifts to characteristic zero.
Since $p\geq P_2'$ the lifted surface and ${\cal Y}_K$
have the same $h^{01}$.
Moreover, these two surfaces have the same $b_1$
by \cite[Lemma 10.2]{ku} and it follows
that $2h^{01}=b_1$ for ${\cal Y}_K$.
In particular, $\Pic^0({\cal Y}_K)$ is reduced.

We choose $P(m)$ to be the maximum of $P_1$ and the $P_2'$'s for
every of the finitely many components of $\cal M$.
Then, every minimal surface of general type with $K^2=m$ over
a field $K$ of characteristic $p\geq P(m)$ corresponds
to a $\Spec K$-valued point of $\cal M$ and we have already
seen that all corresponding surfaces have a reduced
$\Pic^0$.
\qed\medskip


The proof 
does not give an effective bound for $P(m)$.
To find such bounds, a more detailed analysis is needed,
which we now illustrate by determining the
optimal $P(1)$ explicitly.

\begin{Proposition}
 Minimal surfaces of general type with $K^2=1$ have a reduced
 $\Pic^0$ over fields of characteristic $p\geq7$.
 There do exist minimal surfaces of general type with
 $K^2=1$ and non-reduced $\Pic^0$ over fields of 
 characteristic $5$.
\end{Proposition}

\prf
By \cite[Proposition 1.1]{godeaux}, such surfaces fulfill
$1\leq\chi(\OO_X)\leq3$, $p_g\leq2$, $b_1=0$ and $h^{01}\leq1$.
Hence if $\chi(\OO_X)=3$ we necessarily have $p_g=2$ and
$h^{01}=0$ and in particular the $\Pic^0$ of such a surface
is reduced.

In case $\chi(\OO_X)=1$ the $\Pic^0$ is reduced in
characteristic $p\geq7$ by \cite[Corollary 2.6]{godeaux},
which is one of the main results of this article.
The first example of such a surface with non-reduced $\Pic^0$ in
characteristic $5$ is due to Miranda \cite{miranda}, cf.
also \cite[Section 5]{godeaux}.

If $\chi(\OO_X)=2$ we either have
$p_g=1$ and $h^{01}=0$, and such a surface has a reduced $\Pic^0$, 
or  $p_g=2$ and $h^{01}=1$,
in which case the surface has a non-reduced $\Pic^0$,
since $b_1=0$.
However, in this latter case there exists a 
$\mu_p$-, or an $\alpha_p$-torsor above $X$ (depending on whether
Frobenius acts bijectively or trivially on $H^1(\OO_X)$),
and arguing as in the proof of \cite[Theorem 2.4]{godeaux}
we find that such surfaces can only exist in characteristic $2$.
\qed

\end{document}